\documentclass[11pt]{article}
\usepackage{latexsym,amsmath,amssymb}
\date{}
\usepackage[margin=2cm]{geometry}

\numberwithin{equation}{section}
\usepackage{amsthm,bbm}
\usepackage{paralist}
\usepackage[usenames]{color}\usepackage{graphicx}
\definecolor{brown}{cmyk}{0, 0.72, 1, 0.45}
\definecolor{grey}{gray}{0.5}

\allowdisplaybreaks

\newtheorem{maintheorem}{Theorem}

\newtheorem*{conjecture*}{Conjecture}
\newtheorem{theorem}{Theorem}[section]
\newtheorem*{theorem*}{Theorem}
\newtheorem{lemma}[theorem]{Lemma}

\newcommand{\ignore}[1]{}
\newcommand{\one}{\mathbbm{1}}
\def\Pr{\mathbb P}
\def\E{\mathbb E}
\renewcommand{\epsilon}{\varepsilon}
\renewcommand{\a}{\alpha}
\renewcommand{\b}{\beta}
\newcommand{\g}{\gamma}
\renewcommand{\d}{\delta}
\renewcommand{\k}{\kappa}
\newcommand{\m}{\mu}

\def\cE{{\cal E}}

\def\D{\Delta}

\def\L{\Lambda}

\def\cF{{\cal F}}

\def\cG{{\cal G}}
\def\cI{{\cal I}}

\newcommand{\ha}{{\hat{\alpha}}}
\newcommand{\hb}{{\hat{\beta}}}
\newcommand{\hg}{{\hat{\gamma}}}

\newcommand{\hk}{{\hat{\kappa}}}
\newcommand{\hm}{{\hat{\mu}}}

\begin{document}
\title{Random greedy triangle-packing \\ beyond the 7/4 barrier}

\author{ {Tom Bohman\thanks{Department of Mathematical Sciences, Carnegie Mellon
University, Pittsburgh, PA 15213. Email address: {\tt tbohman@math.cmu.edu}.
Research supported in part by NSF grant DMS-1001638.}} \qquad
{Alan Frieze\thanks{Department of Mathematical Sciences, Carnegie Mellon
University, Pittsburgh, PA 15213. Email address: {\tt alan@random.math.cmu.edu}.
Research supported in part by NSF grant DMS-0721878.}} \qquad
{Eyal Lubetzky\thanks{Theory Group of Microsoft Research, One Microsoft Way, Redmond, WA 98052. Email address:
{\tt eyal@microsoft.com}.}}}

\maketitle
\begin{abstract}
The random greedy algorithm for constructing a large partial Steiner-Triple-System is defined as follows.
Begin with a complete graph on $n$ vertices and proceed to
remove the edges of triangles one at a time, where
each triangle removed is chosen uniformly at random
out of all remaining triangles.
This stochastic process terminates once it arrives at a triangle-free graph, and a longstanding open problem is to estimate the final number of edges, or equivalently the time it takes the process to conclude. The intuition that the edge distribution is roughly uniform at all times led to a folklore conjecture that the final number of edges is $n^{3/2+o(1)}$ with high probability, whereas the best known upper bound is $n^{7/4+o(1)}$. It is no coincidence that various methods break precisely at the exponent 7/4 as it corresponds to the inherent barrier where co-degrees become comparable to the variations in their values that arose earlier in the process.

In this work we significantly improve upon the previous bounds by establishing that w.h.p.\ the number of
edges in the final graph is at most $ n^{ 5/3+o(1)} $. Our approach relies
on a system of martingales used to control key graph parameters, where the crucial new idea is to harness the self-correcting nature of the process in order to control these parameters well beyond the point where their early variation matches the order of their expectation.
\end{abstract}

\section{Introduction}

We consider the random greedy algorithm for triangle-packing.
This stochastic graph process begins with the
graph $ G(0) $, set to be the
complete graph on vertex set $[n]$, then proceeds to repeatedly remove the edges of randomly
chosen triangles (i.e.\ copies of
$K_3$) from the graph. Namely, letting $G(i)$ denote the graph that remains after $i$ triangles have been removed,
the $(i+1)$-th triangle removed is chosen uniformly at random from the set of all triangles in $ G(i)$.
The process terminates at a triangle-free graph
$G(M)$. In this work we study the random variable $M$, the number of triangles removed before
obtaining a triangle-free graph (or equivalently, the number of edges in the final triangle-free graph, which is
$ \binom{n}{2} - 3M $).

This process and its variations play an important role in the
history of combinatorics.
R\"odl \cite{rodl} proved the Erd\H{o}s-Hanani conjecture --- which posits the existence of large partial Steiner systems,
collections of $t$-sets with the property that no $k$-set is a subset of more than one set in the collection ---
in the early 1980's by way of a randomized
construction that is now known as the R\"odl nibble.  This construction is a semi-random variation on the
random greedy packing process defined above.  It is semi-random in the sense that the desired object is
constructed in a sequence of substantial pieces, where the proof of the existence of each piece is
an application of the probabilistic method.
Such semi-random constructions have been successfully applied to establish
various key results in combinatorics (see~\cite{AKSz} for an early application of
this approach and~\cite{AS}~and~\cite{MR} for further details). We note in passing that semi-random
techniques have been refined to show the existence of partial Steiner systems that are nearly
as large as allowed by the simple volume bound, see \cite{KR}~and~\cite{van}.
In particular, Alon, Kim and Spencer~\cite{AKS} used such techniques
to prove the existence of a set of
edge-disjoint triangles on $n$ vertices that covers all but $ O( n^{3/2} \log^{3/2} n) $ edges
of the complete graph.

Despite the success of the R\"odl nibble, the limiting behavior of the random greedy
packing process itself remains unknown, even in the special case of triangle
packing
considered here.  Recall that $G(i)$ is the
graph remaining after $i$ triangles have been removed.  Let $ E(i) $ be the edge
set of $ G(i) $.  Note that $ |E(i)| = \binom{n}{2} - 3i $ and that
$ E(M) $ is the number of edges in the triangle-free graph produced by the
process.   It is widely believed that the graph produced by
the random greedy triangle-packing
process behaves like the Erd\H{o}s-R\'enyi random graph with the same edge density, hence the
process should end once its number of remaining edges becomes comparable to the number of
triangles in the corresponding Erd\H{o}s-R\'enyi random graph (i.e., once $|E(M)|$ matches
the order of $(|E(M)|/\binom{n}2)^3 \binom{n}3$). Throughout the paper an event is said to hold \emph{with high probability} (w.h.p.) if its probability tends to $1$ as $n\to\infty$.
\begin{conjecture*}[Folklore]
\label{conj:kahuna}
With high probability $ |E(M)|= n^{ 3/2 + o(1)} $.
\end{conjecture*}
\noindent
Joel Spencer has offered \$200 for a resolution of this question. It was shown by Spencer~\cite{joel} in 1995, and independently by R\"odl and Thoma~\cite{rodl} in 1996,
that $ |E(M)| = o(n^2) $ w.h.p. Grable~\cite{grable} improved this bound to $ |E(M)| \le n^{11/6 + o(1)} $
via an adaptation of the R\"{o}dl nibble method and further sketched how similar arguments using
more delicate calculations should extend this to a bound of $n^{7/4+o(1)}$ w.h.p.
Wormald~\cite{nick2} later demonstrated how the differential equation method can also give nontrivial bounds for this problem (as well as generalizations of it), and namely that $|E(M)| \leq n^{2-\frac{1}{57}+o(1)}$. Finally, in a companion paper~\cite{BFL} that introduced a differential equation approach to this process exploiting its self-correction nature, the foundations of the present work, the authors gave a short proof that $ |E(M)| = O (n^{7/4+o(1)}) $ w.h.p.

It is important to note that the point at which there are roughly $n^{7/4}$ remaining edges is a natural barrier in the analysis
of this process.  To illustrate this, notice that if the $(i+1)$-st triangle taken is $abc$ then the change in the number of triangles in the graph once $abc$ is removed is simply
$ - |N_{a,b}(i)| - |N_{a,c}(i)| - |N_{b,c}(i)| +2$, where $ N_{u,v}(i)$ denotes the common neighborhood of the vertices $u,v \in [n]$ in the graph $G(i)$.
Hence, a natural prerequisite to analyzing this process
is the understanding of the co-degrees $|N_{u,v}|$ for all $u,v$. Suppose for the sake of this
discussion that early in the evolution of the process $G(i)$ closely
resembles the random graph with the same number of edges; that is,
suppose $G(i)$ is roughly the same as $ G_{n,p} $ where $ p = p(i) = 1 - 3i/\binom{n}2$.
If this is the case when $p$ is  close to $1/2$ (i.e.\ $i$ is nearly $ \binom{n}{2}/6$)
then we expect the $|N_{u,v}|$'s  to be close to
$ n/4$ with variations as large as $ \sqrt{n} $.  If these variations in co-degrees persist to the
point where $ p = n^{-1/4} $ (that is, $i$ roughly $ \binom{n}{2}/3 - n^{7/4}$), where we expect
the $|N_{u,v}|$'s themselves to be roughly $ n^{1/2}$, then these variations would be as
large as their average value. Once this happens all control over co-degrees is lost, e.g.\ one could have all co-degrees 0 with non-vanishing probability, or half of the
co-degrees 0 and the other half around $ n^{1/2}$, etc.  In any case, if the variations in $ |N_{u,v}| $ that develop early in the process are not somehow dealt with, one would expect the
analysis to break down once $n^{7/4}$ edges remain.
Perhaps this is the reason that Wormald~\cite{nick2}, who also treated this process
with the differential equation method, stated that ``some non-trivial modification would be required
to equal or better Grable's result.''

In this work we exploit the self-correcting nature of the process in a system of carefully constructed martingales which allows us to tighten the control over key graph properties over time and overcome
the variations in their values that arise early in the process. Our main result is an upper bound on $ |E(M)| $ that is significantly
better than $n^{7/4}$.
\begin{maintheorem}
\label{thm:partial}
Consider the random greedy algorithm for triangle-packing
on $n$ vertices. Let $M$ be the number of steps it takes the algorithm to terminate and
let $E(M)$ be the edges of the resulting triangle-free graph. Then with high probability
$|E(M)| = O\big( n^{ 5/3} \log^{4} n\big)$.
\end{maintheorem}
A key feature of our proof of Theorem~\ref{thm:partial} is an estimate for $ |N_{u,v}| $ in which
the variation {\em decreases} as the process evolves.  We stress that estimates for random graph processes with this property are not commonly obtained by martingale arguments or the differential equation method.

The remainder of the paper is organized as follows. In the next section we discuss our analysis of this
process in more detail, listing the random variables that we track and the estimates on them that
we are able to prove.  The proof of our main result, Theorem~\ref{thm:enchilada}, follows in Section~3.
Theorem~\ref{thm:partial} follows directly from Theorem~\ref{thm:enchilada}.


\section{Self-correcting Estimates}
Let $ (\cF_i) $ be the filtration given by the underlying process.  We note in passing
that our probability space is the set of all maximal sequences of edge-disjoint triangles on
vertex set $ [n] $ with probability measure given by the uniform random choice at
each step.
For $u,v,w \in [n]$ define
$ N_u = N_u(i) = \{ x \in [n] : x u \in E(i) \}$, let $N_{u,v} = N_u \cap N_v$ and let $N_{u,v,w}
= N_u \cap N_v \cap N_w$.
Our main interest is in tracking the number of triangles in $ G(i) $ and the variables
\begin{align*}
Y_{u,v}(i) &= |N_{u,v}(i)| = | \{ x \in [n] : x u \in E(i),\, xv \in E(i) \} |\,.
\end{align*}
In the course of our argument we will also need to consider the variables
\[ Y_u(i) = |N_u(i)| \qquad \mbox{ and } \qquad Y_{u,v,w}(i) = |N_{u,v,w}(i)|\,.
\]
We begin by writing the one-step expected changes in our main variables of interest.
For any random variable $W$ let $ \D W $ be the one-step change $ \D W = W(i+1) - W(i) $.
Let $ Q(i) $ be the number
of triangles in $ G(i) $.  We have
\begin{align}
\E[ \D Y_{u,v} \mid \cF_i ] & = - \sum_{x \in N_{u,v}} \frac{ Y_{u,x} + Y_{v,x} - \one_{\{uv \in E\}}}{ Q} \,,\label{eq-Y-one-step}\\
\E[\D Q \mid \cF_i] &= - \sum_{ xyz \in Q} \frac{ Y_{x,y} + Y_{x,z} + Y_{y,z} - 2}{Q}\,. \label{eq-Q-one-step}
\end{align}
We use these one-step expected changes to relate the random variables to
functions of a continuous `time' (following the approach to the differential
equation method introduced in \cite{r3t}).
We choose the time-scaling $ t = t(i) = i/n^2 $.  Following the
convention established in the
Introduction we set
\begin{equation}
  \label{eq-p-def}
  p = p(i,n) = 1 - \frac{6i}{n^2} = 1 - 6t\,.
\end{equation}
Note that $p$ can now be viewed as either a function of $i$ or the continuous time $t$; we pass
between these interpretations without further comment throughout the paper.
Now, these choices yield the trajectories
$ Y_{u,v}(i) \approx y(t) n $ and $ Q(i) \approx q(t) n $ where we set
\[ y(t) = p^2(t) \qquad \mbox{ and } \qquad q(t) = p^3(t)/6. \]
We can arrive at these equations by either deriving them from
\eqref{eq-Y-one-step}~and~\eqref{eq-Q-one-step} and the assumption that the one-step changes in the trajectory
are equal to the expected one-step change in the corresponding random variable or by
appealing to our $ G_{n,p}$ intuition.
The companion paper \cite{BFL} uses these variables alone to establish a bound
of $O(n^{7/4}\log^{5/4}n)$ on the number of edges that survive to the conclusion of the
algorithm.

In order to achieve better precision, we introduce additional variables
with the central goal of establishing an estimate for $ Y_{u,v}$ with variation that
decreases as the process evolves.  (For applications of the differential equation
method that exploit this kind of `self-correcting' phenomenon, see \cite{mike} and \cite{nick}.)
We would like to add random variables to our collection that will give us better control
on the expression in the numerator of \eqref{eq-Y-one-step}, the one-step expected change in $
Y_{u,v} $.  To this end we take
a closer look at this expression.  We have
\begin{equation}
  \label{eq-D-Yuv-PRS}
  \sum_{x \in N_{u,v}} (Y_{u,x} + Y_{v,x}
- \one_{\{uv \in E\}})= R_{u,v} + R_{v,u} + Y_{u,v}\one_{\{uv\in E\}}
\end{equation}
where
\begin{align*}
 R_{u,v} &= \Big|\left\{ (x,y): xy \in E ~,~ x \in N_{u,v} ~,~
 y \in N_u ~\mbox{ and }~ y \neq v \right\}\Big|\,.
\end{align*}
(Notice that $R_{u,v}$ counts ordered pairs, thus edges in $N_{u,v}$ are counted twice.)
We expect to have $ R_{u,v}  \approx  p\, Y_u \, Y_{u,v} $ which in turn suggests
that $ \sum_{x \in N_{u,v}} (Y_{u,x} + Y_{v,x} - \one_{\{uv \in E\}}) \approx p\, Y_{u,v}( Y_u + Y_v)$.
This expression is in a form that should provide self-correction
in our estimate for $ Y_{u,v} $. Indeed, if $Y_{u,v}$ is large compared to its average
then so will this term be and so (as this term is negated in the
expected one-step change) $Y_{u,v}$ will have a drift back toward its mean.
This approximation emphasizes the need to control vertex degrees: turning to $Y_u$ we have
\[  \E \left[ \D Y_u \mid \cF_i \right] = - \frac{1}{ Q} \sum_{x \in Y_u} Y_{u,x} = - \frac{ 2 T_u }{ Q} \quad\mbox{ where }\quad
T_u =\Big|\left\{ xy \in E(i) : x,y\in N_u \right\}\Big|\,.
\]
The variable $T_u$ again lets us bypass the accumulation of worst case individual errors in the summation of $Y_{u,v}$ variables.
We expect to have $T_u \approx p \binom{ Y_u}{2}$.
Finally, control over triple-degrees $Y_{u,v,w}$ is further needed for our concentration arguments to hold.
We thus arrive at the ensemble of variables $Q, Y_{u,v}, R_{u,v}, Y_u, T_u $
and $ Y_{u,v,w} $ for all $u,v,w\in V_G$. 
The following theorem establishes concentration for this ensemble (throughout the paper $A=B\pm C$ is short for $A\in [B-C,B+C]$).

\begin{theorem}
\label{thm:enchilada}
Set
\[ \g = \frac12 ~,\qquad \hg = \hg(n) = \gamma - \frac{6}{ \log n } \qquad\mbox{ and } \qquad  \Phi = \Phi(p,n) = e^{1-p}\log^2 n\,.\]
Then there exist absolute constants $\a, \b ,\k, \m > 0$
such that with high probability
\begin{align}
Q & = n^3p^3/6  \pm \alpha^2 n^2 p^{2\hg -1} \Phi^2 \label{Qbound}\\
Y_{u,v} & = n p^2  \pm \a n^{1/2} p^{\hg} \Phi  \label{Yuvbound} \\
R_{u,v} & = p Y_u Y_{u,v} \pm  \b n^{3/2} p^{2 + \hg } \Phi  \label{Ruvbound} \\
Y_u & = np  \pm  \kappa n^{1/2}p^{\hg-1} \Phi \label{Yubound}\\
T_u&= pY_u^2/2  \pm  \m n^{3/2}p^{1+\hg} \Phi\label{Tubound}\\
Y_{u,v,w} & = n p^3 \pm  2 \sqrt{n p^3 \log^5 n}\label{Yuvwbound}
\end{align}
for all $ u,v,w $ and as long as $p \geq p^\star = \big(6\alpha^2 e^{2}\, \frac{\log^{10} n}n\big)^{1/(4-2\hg)}$.
\end{theorem}
To deduce Theorem~\ref{thm:partial}, observe that $p^\star$, defined as the smallest $p$ for which the theorem holds, satisfies
\begin{equation}
  \label{eq-p-star}
  p^\star=O\left(n^{-1/3}\log^{10/3}n\right)
\end{equation} since $4-2\hg = 3+O(1/\log n)$. In particular, $p=p^\star$ satisfies $n^3p^3/6 > \alpha^2 n^2 p^{2\hg-1}\Phi^2$ since we have $\Phi = e^{ (1-p)} \log^2 n \leq e \log^2 n$. It thus follows that $Q>0$ w.h.p.\ due to~\eqref{Qbound} and it remains to recover the number of edges corresponding to $p^\star$.
Recalling~\eqref{eq-p-def} we have
\[ |E(i)| = \binom{n}2 - 3i = \binom{n}2 - \frac12 (1-p^\star)n^2 = \frac12 \left(n^2 p^\star - n\right)\]
and the desired result follows from~\eqref{eq-p-star} with room to spare in the power of the logarithmic factor.

We prove Theorem~\ref{thm:enchilada} in the following section by applying martingale arguments
to random variables that track the
differences between the random variables we are interested in and the
variables they should follow.
Note that we establish some form of self-correction for every variable in
this ensemble, with the notable exception of $Y_u$.

The authors suspect that the methods introduced in this paper can be further developed
to achieve better high probability upper bounds on $ |E(M)| $.  This might be achieved
by expanding the ensemble of random variables (perhaps using ensembles of
generalized extension-counting variables, which is the
approach taken in the recent analysis of the $H$-free process \cite{BK}).
However, it seems that a nontrivial modification would be needed to prove the
conjectured bound $ |E(M)| = n^{3/2 + o(1)} $.

For notational convenience we set
\[ \L = \frac{1}{\log^2 n}   \,. \]
Note that while Thereom~\ref{thm:enchilada} applies, estimates \eqref{Yuvbound}--\eqref{Tubound}
and \eqref{Qbound} can each be written as a main term times $ (1 + O( 1/ \log^2 n)) = (1 + O(\L)) $.
Throughout the paper we will use a convention whereby all Greek letters are universal constants.
We do not replace any of the constants (including the pivotal $ \g $) with their
actual values.
This is done in the interest of understanding the role these
constants play in the calculations; it turns out that these constants must be
balanced in a fairly delicate way.  We observe that these constants can take the actual values
\begin{align*}
\a &= 1\,,    &\b & = \tfrac{1}{2}\,, 
&\d& =\tfrac{1}{3}\,,  &\k &=\tfrac{1}{4}\,,  &\m &= \tfrac{1}{4} \,.
\end{align*}
The key conditions these constants must satisfy are \eqref{eq:keyR},
\eqref{eq:keyYuv}, \eqref{eq:keyT} and \eqref{eq:keyY}. 


\section{Proof of Theorem~\ref{thm:enchilada}}


Define $p^\star$ as in~\eqref{eq-p-star} and let $i^\star = \frac16(1-p^\star)n^2$ be the analogous round. Let $ \cG_i $ be the event
that all estimates in Theorem~\ref{thm:enchilada} hold for the first $i$ steps of the process.

For each variable and each bound (upper and lower) in Theorem~\ref{thm:enchilada}
we define a {\em critical interval}.
This interval has one end at the bound we are trying to maintain and the other end slightly closer
to the expected trajectory of the random variable.  If one of the estimates of Theorem~\ref{thm:enchilada} is violated
then the corresponding random variable `crosses' a critical interval.  We bound the probability of each
such event using a martingale argument, introducing a separate supermartingale for each
variable and bound of interest and for each step in which the variable could {\em enter} the critical
interval.  Theorem~\ref{thm:enchilada} then follows from the union bound
(note that the number of supermartingales we consider is bounded by a polynomial in $n$).
We restrict our attention to these critical intervals because the expected one-step
changes in our random variables each have a `drift' term that pushes a wayward
variable back toward the expected trajectory.  By only considering the critical intervals
we make full use of these terms: This is the mechanism we are using to establish
self-correcting estimates.  For an application of this idea in a setting with fewer
variables, see \cite{BFL}.

Let the stopping time $\tau$ be the minimum of $i^\star$ and
the smallest index $i$ such that
$ \cG_i $ does not hold.    Consider an event $ \cE $ of the form $ X(i) \le x(t) $ for all
$ i \le i^\star $ where we assume that $X(i)$ is a random variable and $x(t)$ is not.
Note that every bound in  \eqref{Qbound}--\eqref{Yuvwbound} can be written in this form; that is, the event
$ \{ \tau = i^\star \} $ can be written
\[ \{ \tau = i^\star \} = \bigcap_{\ell \in \cI} \cE_\ell \]
where $ | \cI| $ is polynomial in $n$ and each event $ \cE_\ell $ is of the
form $ X(i) \le x(t) $ for all
$ i \le i^\star $.  For each such event $\cE$ we introduce
a critical interval of the form $ I_\cE(t) =( x(t) - w(t), x(t)) $ where $ w(t) = o(x(t)) $.
Consider a {\em fixed step} $ i_0 $, which we view as a step at which $X(i)$ might enter the
critical
interval $I_\cE$.  Set $ t_0 = i_0/n^2$.  Define the stopping
time $ \tau_{\cE,i_0} $ to be the minimum of $ \max\{i_0,\tau\} $ and
the smallest $ i \ge i_0 $ such that $ X(i) \not\in I_\cE $.  Note that if $ X(i_0) \not\in I_\cE( t_0) $ then
we have $ \tau_{\cE,i_0} = i_0 $.  Thus this stopping time is (formally) well-defined on the
full probability space (n.b.\ we only make use of this stopping time when $ X(i_0) $ is in the critical interval
and $ X(i_0-1)$ is not).
We now establish a bound $B(i)$ on the one-step change in $X(i)$
conditioned on $ \cG_i$. This bound is
far less than the width  $ w(t) $ of the critical interval.  Given a particular event $ \cE$ and starting step
$ i_0$, we work with
the sequence of random variables
$$ Z_{\cE, i_0}(i) = Z(i) = \begin{cases} X(i) - x(t)  & \mbox{if } i_0 \le i  \le \tau_{\cE,i_0} \\
Z(i-1) & \mbox{otherwise.} \end{cases}
$$
Note that if $ X(i_0 -1) $ is not
in the critical interval (and $i_0 <\tau$) then
$ Z(i_0) < - w(t_0) + B(i_0) $.  Therefore, in the event $ \cE^c $ there are
steps $ i_0 < j \le i^\star $ such that $ Z(i_0) < - w(t_0) + B(i_0) $ and $ Z(j) \ge 0 $.
However, our stopping time $\tau$ stops all of these sequences as soon as any of our conditions
\eqref{Qbound}--\eqref{Yuvwbound} are violated.  So, we have
\[
\begin{split}
\{ \tau < i^\star \} & \subseteq \bigcup_{ \ell \in \cI} \bigcup_{ 1 \le i_0 < i^\star }
\{ Z(\tau) - Z( i_0) > w( t_0) - B(i_0) \} \\
& \hskip1cm = \bigcup_{ \ell \in \cI} \,\bigcup_{ 1 \le i_0 < i^\star }
\{ Z(i^\star) - Z( i_0) > w( t_0) - B(i_0) \}\,.
\end{split}
\]
It remains to bound the probability of each event in this union.  This is done for each of the bounds \eqref{Qbound}--\eqref{Yuvwbound}
in turn in Sections~\ref{sec:R}~--~\ref{sec:Q}. In order to bound the probability of these events we will apply the following inequality due to Freedman~\cite{Free}, which was originally stated for martingales yet its proof extends essentially unmodified to supermartingales.
\begin{theorem}[\cite{Free}, Thm~1.6] \label{thm:concentrate}
Let $(S_0,S_1,\ldots)$ be a supermartingale w.r.t.\ a filter $(\cF_i)$.
Suppose that $S_{i+1}-S_i\leq B$ for all $i$, and write $V_t=\sum_{i=0}^{t-1} \E\left[(S_{i+1} - S_{i})^2 \mid \cF_{i}\right]$. Then
for any $s,v>0$
\[\Pr\big(\{S_t\geq S_0+s,\,V_t\leq v\}\mbox{ for some $t$}\big)\leq \exp\bigg\{-\frac{s^2}{2(v+Bs)}\bigg\}\,.\]
\end{theorem}
Our applications of this inequality will each have
two parts: a careful calculation that establishes a martingale condition and
a coarser argument that provides bounds on both
the one-step changes and the second moment of the one-step changes of these variables.
We emphasize that our carefully chosen stopping times
allow us to assume that the event $ \cG_i$ holds throughout
these calculations.  This is henceforth assumed without further comment.

\subsection{Edges between a co-neighborhood and a neighborhood ($R_{u,v}$)}
\label{sec:R}


We begin with an analysis of the one-step expected
change. There are 7 types of triangles that contribute to $\E[\D R_{u,v}\mid \cF_i]$.
\begin{enumerate}[(1)]
\item Triangles $vxy$ where $x\in N_{u,v}$ and $y\in N_v \setminus N_{u,v}$ and $y \neq u$.
For $x \in N_{u,v}$ there are $Y_{v,x}-Y_{u,v,x} - \one_{\{uv \in E\}}$ such triangles
and selection of one of these triangles moves $x$ from $N_{u,v}$ to $N_u \setminus N_{u,v}$ and thereby decreases $R_{u,v}$ by $Y_{u,x}- \one_{\{uv \in E\}}$. This results in a contribution of
\begin{align}
- \sum_{ x \in N_{u,v}} &\frac{ Y_{v,x} - Y_{u,v,x}  - \one_{\{uv \in E\}}}{Q} \left( Y_{u,x} - \one_{\{uv \in E\}}\right)
 = - \sum_{ x \in N_{u,v}} \frac{ Y_{v,x} - Y_{u,v,x} }{Q} Y_{u,x} + O \left( \frac{p}{n} \right)\,,\label{eq-R-type-1}
\end{align}
where in the last equality we absorbed the indicator variables into an $O(p/n)$ term based on the estimates of Theorem~\ref{thm:enchilada}
up to this point.
\item Triangles $uxy$ where $x\in N_{u,v}$ and $y\in N_u\setminus N_{u,v}$ with $y\neq v$. (There are
are $Y_{u,x}-Y_{u,v,x}-\one_{\{uv\in E\}}$ such vertices $y$ for each $x\in N_{u,v}$ and there are $Y_{u,v,y}$ such vertices $x$ for each $y\in N_u\setminus ( N_{u,v}\cup\{v\})$.)
The selection of
such a triangle removes $y$ from $N_u$ while moving $x$ from $N_{u,v}$ to $N_v \setminus N_{u,v}$.
The effect of $y$ is a decrease of $Y_{u,v,y}$ in $R_{u,v}$ while the effect of $x$ is an additional decrease of $Y_{u,x} + Y_{u,v,x} - \one_{\{uv \in E\}} - 1$
(the vertex $y\in N_{u,x}$ was already counted). The overall contribution is
\begin{align}
- \sum_{ x \in N_{u,v}} \frac{ Y_{u,x} - Y_{u,v,x} }{Q}(Y_{u,x}+Y_{u,v,x})
- \sum_{ y \in N_u \setminus N_{u,v}} \frac{ Y_{u,v,y}^2}{Q} + O \left( \frac{p}{n} \right)\,,\label{eq-R-type-2}
\end{align}
where the indicator terms were again absorbed into the $O(p/n)$ term.
\item Triangles $vxy$ where $x,y\in N_{u,v}$. Choosing such a triangle moves $x,y$ from $N_{u,v}$ to $N_u\setminus N_{u,v}$ and thus
decreases
$R_{u,v}$ by $(Y_{u,x}- \one_{\{uv \in E\}}) +(Y_{u,y}-\one_{\{uv \in E\}})$. This contributes
\begin{align}
-\sum_{\substack{x,y\in N_{u,v}\\ xy \in E}}\frac{1}{Q}(Y_{u,x}+Y_{u,y}) + O \left( \frac{p^2}{n} \right) =
- \sum_{x \in N_{u,v}} \frac{ Y_{u,v,x}}{Q}  Y_{u,x}  + O \left( \frac{p^2}{n} \right)\,.\label{eq-R-type-3}
\end{align}
\item Triangles $uxy$ where $x,y\in N_{u,v}$. Selection of such a triangle
moves $x,y$ from $N_{u,v}$ to $N_v\setminus N_{u,v}$ and so decreases
$R_{u,v}$ by $(Y_{u,x}+Y_{u,v,x} - \one_{\{uv \in E\}})  +(Y_{u,y}+Y_{u,v,y}- \one_{\{uv \in E\}})-2$, translating to
\begin{align}
 -\!\!\!\sum_{\substack{x,y\in N_{u,v}\\ xy \in E}}\frac{1}{Q}(Y_{u,x}+Y_{u,v,x}+Y_{u,y}+Y_{u,v,y}) + O\left( \frac{p^2}{n} \right)
= -\sum_{x \in N_{u,v}} \frac{ Y_{u,v,x}}{Q} ( Y_{u,x} + Y_{u,v,x}) + O \left( \frac{p^2}{n} \right).\label{eq-R-type-4}\end{align}
\item Triangles $uxy$ where $xy$
is in the set of edges induced by
$N_u\setminus N_{u,v}$. Each such triangle decreases
$R_{u,v}$ by $Y_{u,v,x}+Y_{u,v,y}$.
This results in a contribution of
\begin{align}-\!\!\!\sum_{\substack{x,y\in N_u\setminus N_{u,v}\\ xy\in E}}\frac{1}{Q}(Y_{u,v,x}+Y_{u,v,y})=-\sum_{x\in N_u\setminus N_{u,v}}\frac{Y_{u,x}-Y_{u,v,x}}{Q}
Y_{u,v,x}\,.\label{eq-R-type-5}\end{align}
\item Triangles $uvx$ for $x\in N_{u,v}$ contribute
 \begin{align}
 -\one_{\{uv\in E\}}\sum_{x \in N_{u,v}}\frac{Y_{u,x} + Y_{u,v,x} -1}Q = -\frac{R_{u,v}\one_{\{uv\in E\}}}{Q} = O \left( \frac{p}{n} \right)\,.\label{eq-R-type-6}
  \end{align}
\item Triangles $xyw$ where $x\in N_{u,v}$ and $y\in N_u$ while $w\neq u$.  We note
in passing
that this is the only type of triangle whose selection impacts $R_{u,v}$ while changing neither
$Y_u - Y_{u,v} $ nor $Y_{u,v}$. For a fixed $xy\in E$ with $x\in N_{u,v}$ and $y\in N_u$
there are $Y_{x,y}-1$ such triangles, and each would decrease $R_{u,v}$ by either 1 or 2.
%
Merely applying our bounds on each $ Y_{x,y} $ term individually would produce
an undesirable error.  Instead, we will sum over $ x \in N_{u,v} $ and use our error
bounds on $ R_{x,u} $.  This should give a better estimate by aggregating multiple $ Y_{x,y} $ terms for
better cumulative error bounds, and indeed this proves to
be a crucial choice.
A triangle $xyw$ in this category will decrease $R_{u,v}$ by 1 if $uw\notin E$ and by $2$ if $uw\in E$.
Recall that $R_{x,u}$ counts the number of edges between $N_x\setminus N_u$ and $N_{x,u}$ plus twice the
number of edges within $N_{x,u}$. Hence, the contribution in this case is precisely $(-1/Q)$ times
\begin{align*}
 \Big(\sum_{x \in N_{u,v}} R_{x,u}\Big) -  R_{v,u}\one_{\{uv \in E\}} &=
  -\sum_{x \in N_{u,v}} \left[p Y_x Y_{u,x} \pm \beta n^{3/2} p^{2 +\hg}\Phi \right] + O \left(n^2 p^4 \right) \\
&=\sum_{x \in N_{u,v}} \left[ p Y_x Y_{u,x}\right] \pm \left(1 + O \left( \L \right) \right)\b n^{5/2}p^{4+\hg} \Phi\,,
\end{align*}
where in the last equality we absorbed the $O(n^2 p^4)$ error term using the fact that $\Lambda = \log^{-2}n$ is $O(n^{1/2} p^\hg)$ (with room to spare).
Plugging in the fact that $Y_x = np \pm \kappa n^{1/2} p^{\hg-1}\Phi$ and using the identity $R_{u,v} = \sum_{x\in N_{u,v}} (Y_{u,x} - \one_{\{u v \in E\}})$,
we can conclude that the contribution in this case is
\begin{align}
-\frac{R_{u,v} np^2}Q \pm  \frac1{Q}\left[ (1 + O( \L ))\kappa n^{5/2} p^{4+\hg} \Phi + (1 + O(\L))\beta n^{5/2} p^{4+\hg} \Phi\right]\,,\label{eq-R-type-7}
\end{align}
where we absorbed all indicator variables (at most $\frac{np^2}Q Y_{u,v} = O(p/n)$) into the final error term.
\end{enumerate}
Now that we have an expression (albeit in 7 parts) for the expected change in $R_{u,v}$,
we are ready to consider the expected change in $R_{u,v}$ relative to its expected trajectory.
Define
\[ X =  R_{u,v} - p Y_u  Y_{u,v} \]
and consider $ \E[ \D X \mid \cF_i] $.  Write
\begin{align}
\D \left[  p Y_u Y_{u,v} \right]
&= p(i+1)  \D \left[ Y_u Y_{u,v} \right]  +
\D p   Y_u(i) Y_{u,v}(i) \,.
\label{eq:changedivide}
\end{align}
We will see that the expected change in $R_{u,v}$ due to triangles of types 1--6 will balance off with the first term
in~\eqref{eq:changedivide}, while the expected change due to triangles of type 7 will be balanced
by the second term.

Collecting \eqref{eq-R-type-1}--\eqref{eq-R-type-6}, the total contribution to $ \E[ \D R_{u,v}\mid \cF_i]$ from triangles of types 1--6 equals
\begin{align}\label{eq1}
 - \frac{1}{Q} \sum_{x \in N_{u,v}} Y_{u,x} \left( Y_{u,x}+Y_{v,x} \right)
 - \frac{1}{Q} \sum_{ x \in N_u }  Y_{u,x} Y_{u,v,x} + O \left( \frac{p}{n} \right).
\end{align}
Furthermore, we can analyze the change in $\Xi = Y_{u,v} Y_u$ by considering the following 3 cases: \begin{compactenum}[(i)]
  \item \label{it-Xi-change-1} Selecting a triangle $uxy$ for $x,y\in N_u\setminus N_{u,v}$: here $\D Y_{u,v}=0$ while $\D Y_u= -2$
  and so $\D\Xi = -2Y_{u,v}$.
  \item \label{it-Xi-change-2} Selecting a triangle $vxy$ for $x\in N_{u,v}$ and $y\neq u$: the co-neighborhood loses $x$ and in addition loses $y$ if $y \in N_{u,v}$, while $Y_u$ remains unchanged.
Thus $\D \Xi = -2Y_u$ if $y\in N_{u,v}$ and $\D \Xi = -Y_u$ otherwise.
  \item \label{it-Xi-change-3} Selecting a triangle $uxy$ for $x\in N_{u,v}$: If $y \in N_{u,v}$ then $\D Y_{u,v} = \D Y_u = -2$ and $\D \Xi = -2Y_u -2 Y_{u,v} +4$. Similarly, if $y=v$ then $\D \Xi = -Y_u-Y_{u,v}+1$. Otherwise, $y \in N_{u}\setminus N_{u,v}$
   and $\D Y_{u,v} = -1 $ while $\D Y_u = -2$, hence $\D \Xi = -Y_u - 2Y_{u,v} + 2$.
\end{compactenum}
Altogether, we can obtain the factors of $2$ in Item~\eqref{it-Xi-change-1} and in the case $y\in N_{u,v}$ of Items~\eqref{it-Xi-change-2},\eqref{it-Xi-change-3} automatically by symmetry when summing over $x$ as follows:
\begin{align*}
\E \big[ \D  &[Y_{u,v} Y_u  ] \mid \cF_i\big] =
  - Y_{u,v} \sum_{x \in N_u \setminus N_{u,v}} \frac{ Y_{u,x}-Y_{u,v,x}}{Q} - Y_u  \sum_{ x \in N_{u,v}} \frac{ Y_{v,x}-\one_{\{u v \in E\}}}{Q} \\
& - \bigg[( Y_u + Y_{u,v} + O(1)) \Big(\sum_{x \in N_{u,v}} \frac{ Y_{u,x} }{Q}\Big)  +
\big( (Y_u+2Y_{u,v})-(Y_u + Y_{u,v})+O(1)\big) \Big(\sum_{y\in N_u \setminus N_{u,v}} \frac{Y_{u,v,y}}Q \Big)\bigg]\,.
\end{align*}
All the triple-degree terms cancel out and we can collect all the $O(1)$-terms and rewrite the above as
\begin{align*}
-\frac1{Q} \sum_{ x \in N_{u,v}}
Y_u (Y_{u,x} + Y_{v,x} )
  -\frac1{Q}  \sum_{x \in N_u } Y_{u,v} Y_{u,x}  +  O \left(\frac1n\right)\,.
\end{align*}
Notice that when multiplying the above by $p$ the error term becomes an additive $O(p/n)$ while the main terms are of order $O(p^2)$. As such, the same estimate holds for the result of multiplying the above by $p(i+1)$ (which differs from $p(i)$ by an additive $O(n^{-2})$ error and thus introduces an extra $O((p/n)^2)$ error term). We can now combine this with the change in $R_{u,v}$ given in \eqref{eq1} to get that the contribution to $ \E [ \D X\mid \cF_i] $ from triangles of types~1--6  and the first term in~\eqref{eq:changedivide} is
\begin{align}
& -\frac{1}{Q} \sum_{x \in N_{u,v}}   \left( Y_{u,x}+Y_{v,x} \right)\left(Y_{u,x} - pY_{u}\right)
 -\frac{1}{Q} \sum_{ x \in N_u } Y_{u,x} (Y_{u,v,x}  - p Y_{u,v}) + O \left( \frac{p}{n} \right)\,.\label{eq-DeltaR-types-1-6}
\end{align}
In order to rewrite the last two summations, we need the following straightforward estimate:
\begin{lemma}
\label{lem:jpatrick}
Let $ (x_i)_{i\in I} $ and $ (y_i)_{i \in I}$ such that
$|x _i - x_j| \le \delta_1$ and $|y_i-y_j| < \delta_2$ for all $i,j \in I$.  Then
\[ \bigg| \sum_{i \in I} x_i y_i -  \frac{1}{|I|} \Big( \sum_{i \in I} x_i \Big)
\Big( \sum_{i \in I} y_i \Big) \bigg|  \le  |I| \delta_1 \delta_2 \,.\]
\end{lemma}
The key observation here is that the first and second summations in~\eqref{eq-DeltaR-types-1-6}
feature the random variable $ X $ itself, a fact which our self-correction argument for $X$ hinges on.
Namely, by definition of $R_{u,v}$ we have $\sum_{x\in N_{u,v}}( Y_{u,x} -\one_{\{uv\in E\}}-  p  Y_u ) = R_{u,v} - pY_{u} Y_{u,v} = X$ and
similarly $\sum_{ x \in N_{u}} ( Y_{u,v,x}  - p Y_{u,v} ) = X$.
By the lemma above and our error estimates from Theorem~\ref{thm:enchilada}, the first summation in~\eqref{eq-DeltaR-types-1-6} is equal to
\[ -\frac{(R_{u,v}+R_{v,u})X}{Q\, Y_{u,v}} + O\left(n^{-1} p^{2\hg-1} \Phi^2 \right)
= -\frac{12+O(\Lambda)}{n^2 p}X + O\left(\frac{\Phi^2}n\right)
\,, \]
where the last equality used our $(1+O(\L))$-approximation
for $(R_{u,x}+R_{v,x})$, $Q$ and $Y_{u,v}$. Similarly, the second summation in~\eqref{eq-DeltaR-types-1-6} can be estimated by
\[ -\frac{2 T_u X}{Q\, Y_{u}} + O\left(\frac{np}Q\cdot n^{1/2}p^{\hg}\Phi \cdot  \sqrt{n p^3 \log^5 n} \right)
= -\frac{6+O(\Lambda)}{n^2 p}X + O\left(\frac{\Phi \log^{5/2}n}n\right)
\,. \]
Altogether, the contribution to $ \E[ \D X \mid \cF_i] $ from triangles of types 1--6 and the first term in~\eqref{eq:changedivide} is
\begin{align}
 - &\frac{ 18  + O( \L ) }{ n^2 p }
X  +  O \left( \frac{\log^5 n}{n} \right) \,.\label{eq-D-X-part-1}
\end{align}

We now turn to the triangles of type 7. As we noted above, we balance the term~\eqref{eq-R-type-7}
with the second term in \eqref{eq:changedivide}, i.e.\ the expected change in
$ p Y_{u} Y_{u,v} $ due to the change in $p$ (which deterministically decreases by $6/n^2$).  The total contribution to
$\E[ \D X \mid \cF_i]$ from these terms is
\begin{align}
& - \frac{R_{u,v} np^2}{Q} + \frac{6}{n^2} Y_u Y_{u,v} \pm \frac{1}{Q}
(1 + O( \L ))(\b + \k) n^{5/2} p^{4+\hg}\Phi = - \frac{ 6+O(\L)  }{ n^2 p}X \pm  (6 + O( \L ))(\b + \k) \frac{ p^{1 + \hg}\Phi}{ n^{1/2}}
 \label{eq-D-X-part-2}\,.
\end{align}
The combination of~\eqref{eq-D-X-part-1},\eqref{eq-D-X-part-2} gives
\begin{align}
\label{eq:Ruvmain}
\E \left[ \D X\mid \cF_i \right] = - \frac{ 24 + O( \L)}{ n^2p} X  \pm  ( 6 +O( \L ))
(\b + \k) \frac{ p^{1 + \hg}\Phi}{ n^{1/2}}\,,
\end{align}
where the $O(\Lambda)n^{-1/2} p^{1+\hg}\Phi $ term absorbed the error-term in~\eqref{eq-D-X-part-1}
by the choice of $p$ in~\eqref{eq-p-star}.

We are now ready to establish the concentration of $R_{u,v}$ via a martingale argument. As outlined above, we introduce two critical
intervals for the random variable $X$, corresponding to the upper bound and lower bound on $R_{u,v}$.
These intervals have one endpoint at the bound we are trying to establish and the other somewhat
closer to zero (corresponding to the expected trajectory of $X$).  For the variable $R_{u,v}$ to violate Eq.~\eqref{Ruvbound} it must be that $X$ crosses one of the critical intervals.

Our critical interval for the upper bound on $R_{u,v}$ is
\begin{equation}
\label{eq:Rcrit}
I_R = \left( \hb n^{3/2} p^{2+\hg}\Phi~,~ \beta n^{3/2} p^{2+\hg}\Phi \right) \qquad \mbox{ where } \qquad \hb = \left( 1- \log^{-1} n \right) \beta\,.
\end{equation}
Suppose that $X(i_0)$ enters $I_R$ for the first time at round $i_0$ (within the time range covered by Theorem~\ref{thm:enchilada}) and define the stopping time $\tau_{R}=\min\{ i \geq i_0 : X(i) < \hb n^{3/2} p^{2+\hg}\Phi\}$, i.e.\ the first time beyond $i_0$ at which $X$ exits the interval through its lower endpoint.
We claim that $Z(i \wedge \tau_{R})$ is a supermartingale, where
\[ Z(i) = X(i) -  \beta n^{3/2} p^{2 + \hg} \Phi \qquad \mbox{for $i \geq i_0$}\,.\]
To see this, write $t = i/n^2$ according to which $p = 1-6t$ and $\Phi = e^{6 t} \log^2 n$, and note that for any $\hg >0$ the second derivative of $f(t) = e^{6 t} (1-6t)^{2+\hg}$ is uniformly bounded in $[0,\frac16]$, hence
\begin{equation}
   \label{eq-D-2nd-Z-term}
 f\big(t+n^{-2}\big) = f(t) + \left(6 e^{6 t} (1-6t)^{2+\hg} - 6(2+\hg)e^{6 t}(1-6t)^{1+\hg}\right)n^{-2} + O(n^{-4})\,.
 \end{equation}
This provides an estimate for $\Delta[\beta n^{3/2} p^{2 + \hg} \Phi]$ between $Z(i+1)$ and $Z(i)$. At the same time, for any $X(i)$ satisfying~\eqref{eq:Rcrit} we can plug in the lower bound this gives for $X$ in~\eqref{eq:Ruvmain} and obtain
\begin{align*}
\E[ &\D Z \mid \cF_i\,,\,\tau_{R} > i] \leq
 - \big( 24 + O( \L ) \big) \hb \frac{ p^{1 + \hg} \Phi}{ n^{1/2}}
  +  \big( 6 +O( \L )\big) (\b+\k) \frac{ p^{1 + \hg}\Phi}{ n^{1/2}} \\
&\qquad\qquad\qquad + \left(1+O(n^{-2})\right)\bigg[ 6 (2+\hg) \beta \frac{ p^{1+\hg}\Phi}{ n^{1/2}} - 6 \beta \frac{ p^{2 + \hg} \Phi}{ n^{1/2}}\bigg] \\
&= \big(-4\hb + 3\b+\k+\hg \b+O(\L)\big)\frac{ 6p^{1 + \hg} \Phi}{ n^{1/2}} + \big(-6\beta +O(\L)\big)\frac{ p^{2 + \hg} \Phi}{ n^{1/2}}\,.
\end{align*}
Since $\hb=(1-\log^{-1} n)\beta$ and $\hg=\g - 6\log^{-1}n$ we have
\[  -4\hb + 3\b+\k+\hg \b + O(\L)=  -(1-\g)\beta + \kappa -\frac{2\beta}{\log n} + O(\L) < -(1-\g)\beta + \kappa-\frac{\beta}{\log n}\,,
\]
where in the inequality absorbed the $O(\L)$ term into a single $(\beta/\log n)$-term for large enough $n$.
Altogether, we conclude that if $n$ is large enough then $Z(i\wedge \tau_{R})$ is a supermartingale provided that
\begin{equation}
  \label{eq:keyR-fine}
  (1-\gamma)\beta + \frac{\beta}{\log n} \geq \kappa
\end{equation}
and $\Lambda = o(\beta)$. In particular, we can relax this condition into the requirement that $\Lambda=o(\beta)$ and
\begin{equation}
\label{eq:keyR}
(1-\g)\b \ge \k\,.
\end{equation}

To apply Freedman's inequality we need to obtain bounds on $ |\D Z| $ and $\E[ \left( \D Z \right)^2 \mid \cF_{i}]$.
Recall that $\D Z = \D R_{u,v} - \D[p\, \Xi] - \D[\b n^{3/2}p^{2+\hg}\Phi]$ where $\Xi = Y_u Y_{u,v}$.
By~\eqref{eq:changedivide} and the fact that $\D p = -6/n^2$ we have $\D[p \,\Xi] = p \D \Xi + O(p^3)$, while
 $\D[\b n^{3/2}p^{2+\hg}\Phi] = O(n^{-1/2}p^{1+\hg}\Phi) = o(1)$. Hence, $\D Z$ will be dominated by $\D R_{u,v} - \D[p\, \Xi]$. There are four cases to consider here:
\begin{compactenum}[(i)]
\item Choosing a triangle that includes $u$ or $v$ and some vertex $x \in N_{u,v}$ (triangles of types 1--4,6 in the analysis of $ \E [ \D R_{u,v} \mid \cF_i] $ above):
There are $ O(n^2 p^4) $ such triangles and selecting one of them affects both $ R_{u,v} $ and $ \Xi $.  As we next specify, the principle terms in these changes are identical and so $|\D Z|$ is bounded by the error terms in our approximations. Indeed, going back to the analysis of the triangle types as well that of $\D\Xi$ we recall the various triangle types satisfied:
 \begin{compactitem}
   \item Type 1: $\D R_{u,v} = -Y_{u,x}+O(1)$ vs.\ $\D \Xi = -Y_u$.
   \item Type 2: $\D R_{u,v} = -(Y_{u,x}+Y_{u,v,x}+Y_{u,v,y})+O(1)$ vs.\ $\D \Xi = -Y_u -2Y_{u,v}+ O(1)$.
   \item Type 3: $\D R_{u,v} = -(Y_{u,x}+Y_{u,y})+O(1)$ vs.\ $\D \Xi = -2Y_u$.
   \item Type 4: $\D R_{u,v} = -(Y_{u,x}+Y_{u,v,x}+Y_{u,y}+Y_{u,v,y})+O(1)$ vs.\ $\D \Xi = -2Y_u - 2Y_{u,v}+O(1)$.
   \item Type 6: $\D R_{u,v} = -Y_{u,x}-Y_{u,v,x}+O(1)$ vs.\ $\D \Xi = -Y_u - Y_{u,v}+O(1)$.
 \end{compactitem}
In all of the above cases the main terms of $\D R_{u,v}$ cancel those of $p\D \Xi$ at the cost of an $O( n^{1/2} p^\hg \log^2 n) $-error due to the approximations of $Y_{u,v}$ and $Y_u$ (dominating all other errors).
 \item Choosing a triangle that includes $u$ but no vertex in $ N_{u,v}$ (triangles of type 5 above):
 there are $O(n^2 p^3)$ such triangles and each corresponds to $\D R_{u,v} = -(Y_{u,v,x}+Y_{u,v,y})$ vs.\ $\D \Xi = -2Y_{u,v}$.
 Hence, in this case $|\D Z|= O( \sqrt{np^3 \log^5 n})$ dominated by triple co-degrees (recalling that $\gamma=\frac12$).
\item Choosing triangles that affect the value of $R_{u,v}$ but contain neither $u$ nor $v$ (type~7 triangles):
Each of these $ O(n^3 p^6) $ triangles corresponds to $|\D Z| = O(1)$ as $\D R_{u,v} = -1$ and $\D \Xi = 0$.
\item Choosing any other triangle: as $R_{u,v}$ and $\Xi $ are both unchanged, these triangles can modify $Z$ by at most $O(p^3)$ due to the additive $6/n^2$ change in $p$.
\end{compactenum}
 The $L^\infty$ bound on $\D Z$ is clearly dominated by rounds of the first sort and $|\Delta Z| = O(n^{1/2} p^\hg \log^2 n)$. For an $L^2$ bound notice that the 4 round types contribute $ O(p^{1+2\hg}\log^4 n)$, $O(p^3 \log^5 n)$, $O(p^3)$ and $O(p^6)$ respectively to $\E[ (\D Z)^2\mid\cF_i]$. As $\hg < 1$ we have $p^3 \log^5 n = o(p^{1+2\hg} \log^5 n)$ while the fact that $p\geq p^\star$ (which has order $n^{-1/(4-2\hg)+o(1)}$ as was seen in Eq.~\eqref{eq-p-star}) implies that
 $\frac{\log^{5}n}n =o(p^{1+2\hg})$ since $n p^{1+2\hg} \geq n^{(3-4\g)/(4-2\g)-o(1)} > n^{1/5}$ for large enough $n$. Altogether it follows that $\E[ (\D Z)^2\mid\cF_i] = O(p^{1+2\hg}\log^4 n)$. Clearly the $L^\infty$ and the $L^2$ bounds on $\D Z$ also hold in the conditional space given $\tau_R > i$.

Recall that we are interested in $Z(i)$ starting at time $i_0$, i.e.\ immediately after $X$ enters the critical interval $I_R$. Let $p_0 = p(i_0) = 1 - 6i_0/n^2$ and observe that our bound on $|\D Z|$ guarantees that $0 \leq X(i_0)-\hb n^{3/2}(p_0)^{2\hg}\Phi \leq O(n^{1/2} (p_0)^\hg \log^2 n)$. Hence,
\[ Z(i_0) \leq (\hb-\b) n^{3/2} (p_0)^{2 + \hg} \Phi + O(n^{1/2} (p_0)^\hg \log^2 n) \leq -\tfrac12 \beta n^{3/2} (p_0)^{2 + \hg} \log n \,,\]
where the final factor of $\frac12$ readily cancels the $O(n^{1/2} (p_0)^\hg \log^2 n)$-term for sufficiently large $n$ since $n (p_0)^2 \log^{-1} n \geq n^{1-o(1)}(p^\star)^2 > n^{1/4}$. We now apply Freedman's inequality (Theorem~\ref{thm:concentrate}) to the supermartingale $ S_j = Z\big((i_0 + j) \wedge \tau_R\big) $ while noting that the above analysis implies that
\begin{align*}
&S_0 \leq -\tfrac12 \beta n^{3/2} (p_0)^{2 + \hg} \log n \qquad\,,\qquad
\max_j |S_{j+1}-S_j| = O\left( n^{1/2} (p_0)^\hg \log^2 n \right) \,, \\
&\sum_{j}\E\big[(S_{j+1}-S_j)^2\mid \cF'_j\big] = O\bigg(\log^5 n \sum_{i\ge i_0}  (p(i))^{1+2\hg }\bigg)
= O \left( n^2 (p_0)^{ 2+2 \hg } \log^5 n \right)\,,
\end{align*}
where $\cF'_j = \cF_{i_0+j}$. We deduce that for some fixed constant $c>0$,
\begin{align*}
\Pr\left( \cup_{j}\big\{ S_j \geq 0\big\} \right) &\leq \exp
\left( - c \frac{  n^{3} (p_0)^{4+2\hg} \log^2 n}{ n^2 (p_0)^{2 + 2\hg} \log^5 n + n^2 (p_0)^{2 + 2\hg} \log^3 n  } \right) = \exp \left( - c n (p_0)^{2} \log^{-3} n  \right)\,,
\end{align*}
which is sufficiently small to afford a union bound over all $u,v$ and time steps $i_0$, hence w.h.p.\ $X$ never crosses the critical interval $I_R$ and so $X(i)\leq \b n^{3/2}p^{2+\hg}\Phi$ for all $u,v$ and $i$. The same argument
shows that w.h.p.\ $X(i)\geq -\b n^{3/2}p^{2+\hg}\Phi$ for all $u,v,i$ by considering the lower interval $\big( - \beta n^{3/2} p^{2 + \hg}\Phi\,,\, - \hb n^{3/2} p^{2 + \hg}\Phi  \big)$
and analyzing the variable $ Z(i) = - X(i) - \beta n^{3/2} p^{2 + \hg}\Phi $. This completes the proof of Eq.~\eqref{Ruvbound}.

\subsection{Co-degrees ($ Y_{u,v}$)}

Following the lines of \S\ref{sec:R} we will establish concentration for
\[ X = Y_{u,v} - n p^2\,.\]
As $p$ decreases by $6/n^2$ with each step, Eq.~\eqref{eq-Y-one-step} and~\eqref{eq-D-Yuv-PRS} show that
\begin{align*}
\E[ \D X\mid \cF_i] &= - \sum_{ x \in N_{u,v}} \frac{ Y_{u,x} + Y_{v,x} - \one_{\{uv\in E\}} }{Q} + \frac6{n}\left(2p-6/n^2\right)  = - \frac{ R_{u,v} + R_{v,u} }{ Q} + \frac{ 12p}{n} +  O\bigg( \frac{1}{ n^2p} \bigg) \,,
\end{align*}
where the last error term absorbed the indicators and $O(n^{-3})$ from the first expression. Substituting our estimates~\eqref{Ruvbound}
for $R_{u,v},R_{v,u}$ we get that this is equal to
\[ - \frac1Q \left( p (Y_{u} + Y_v) Y_{u,v} \pm 2 \beta n^{3/2}p^{2+\hg}\Phi\right) +
\frac{12p}{n} + O \bigg( \frac{1}{ n^2p}  \bigg)\,,\]
and using the estimate~\eqref{Yubound} for $Y_u$ and that $Q = (1+O(\L))n^3p^3$ we can conclude that
\begin{align*}
\E[ \D X\mid \cF_i] & = - \frac{ p ( 2np \pm 2\k n^{1/2}p^{\hg-1}\Phi) Y_{u,v}}{ Q } +
\frac{12p}{n}
\pm (12+O(\L))\frac{ \beta p^{\hg-1}\Phi}{ n^{3/2}} +  O \bigg( \frac{1}{n^2p} \bigg)\,.
\end{align*}
Crucially, we did not approximate $Q$ in the first expression with a $(1+O(\L))$ correction factor as this would incur an error that would be too large to handle. Instead, there we apply Eq.~\eqref{Qbound} and the fact that $Y_{u,v} = X + np^2$ to obtain that
\begin{align*}
\frac{ p ( 2np \pm 2\k n^{1/2}p^{\hg-1}\Phi) Y_{u,v}}{ Q } &= \frac{12+O(\L)}{n^2 p}X +
\frac{ 2n^2 p^4 }{\frac16 n^3 p^3 \pm \a^2 n^2 p^{2\hg-1}\Phi^2} \pm
(12+O(\L))\frac{ \k p^{\hg-1}\Phi}{n^{3/2}} \\
&= \frac{12+O(\L)}{n^2 p}X + \frac{12p}n \bigg(1+O\bigg(\frac{p^{2\hg-4}\Phi^2}n\bigg)\bigg) \pm (12+O(\L))\frac{ \k p^{\hg-1}\Phi}{n^{3/2}}\,.
\end{align*}
Combining this with the above estimate for $\E[\D X\mid\cF_i]$, the term $12p/n$ vanishes and we get that
\begin{align*}
\E[ \D X\mid \cF_i] & = -\frac{12+O(\L)}{n^2p}X \pm (12+O(\L))\bigg[ \frac{ \beta p^{\hg-1}\Phi}{ n^{3/2}} +\frac{ \k p^{\hg-1}\Phi}{n^{3/2}}\bigg] +  O \bigg(\frac{p^{2\hg-3}\Phi^2}{n^2}\bigg) \,,
\end{align*}
where the $O\big(1/(n^2 p)\big)$ error term was absorbed into the $O\left(n^{-2}p^{2\hg-3}\Phi^2\right)$-term since $\hg\leq 1$ and so $1/(n^2 p) = o(p^{2\hg-3}\Phi^2 / n^2)$. Furthermore, we claim that one may now omit this latter error-term altogether as it is negligible compared to the error-term of $O(\L)$ in the terms involving $\b,\k$. Indeed, keeping in mind that $\Phi$ and $\L^{-1}$ are each of order $\log^2 n$, we have
\begin{equation} \frac{p^{2\hg-3}\Phi^2}{n^2} = O(\L)\frac{p^{\hg-1}\Phi}{n^{3/2}} \frac{\log^4 n}{n^{1/2}p^{2-\hg}} \leq
O(\L)\frac{p^{\hg-1}\Phi}{n^{3/2}} \frac{\log^4 n}{\sqrt{n(p^\star)^{4-2\hg}}} =
O(\L)\frac{p^{\hg-1}\Phi}{n^{3/2}} \frac{\log^4 n}{\log^5 n} = o\bigg(\L\frac{p^{\hg-1}\Phi}{n^{3/2}}\bigg)\,.
\label{eq-p-gamma-comparison}
\end{equation}
Assume now that $X(i)$ enters the upper critical interval defined by
\[
I_Y = \left( \ha n^{1/2} p^{\hg}\Phi~,~ \a n^{1/2} p^{\hg}\Phi \right) \qquad \mbox{ where } \qquad \ha = \left( 1- \log^{-1} n \right) \a\,.
\]
That is, suppose that $i_0$ is the first round at which $X(i) \geq \ha n^{1/2} p^{\hg} \Phi$ and define the stopping time $\tau_Y = \min\{ i \geq i_0 : X(i) < \ha n^{1/2} p^\hg \Phi\}$. As before, consider
\[ Z(i) = X(i) - \alpha n^{1/2} p^\hg \Phi \,.\]
By the same argument of~\eqref{eq-D-2nd-Z-term} we have
\[ \D\left[\a n^{1/2} p^{\hg} \Phi\right] = \left(6+O(n^{-2})\right)\left( p^{\hg} - \hg p^{\hg-1}\right)\a n^{-3/2} \Phi \]
and combined this with the above upper bound on $\E\left[ \D X\mid \cF_i\right]$ establishes that
\begin{align*}
\E\left[ \D Z \mid \cF_i\,,\,\tau_Y > i\right] & \leq \big[-2\ha + 2\b + 2\k + \hg\a +O(\L)\big]\frac{ 6p^{\hg-1}\Phi }{ n^{3/2}}
+ \big[ -\a + O(\L)\big] \frac{ 6p^\hg \Phi }{ n^{3/2}} \\
&\leq \big[-(2-\g)\a + 2\b + 2\k \big]\frac{ 6p^{\hg-1}\Phi }{ n^{3/2}}
+ \big[ -\a + O(\L)\big] \frac{ 6p^\hg \Phi }{ n^{3/2}} \,,
\end{align*}
where we used the fact that $\g - \hg = 6/\log n$ to absorb both $2(\a-\ha)=2\a/\log n$ and the $O(\L)$-term for large $n$. In particular, $S_j = Z((i_0+j) \wedge \tau_Y)$ is indeed a supermartingale so long as
$\Lambda=o(\alpha)$ and
\begin{equation}
\label{eq:keyYuv}
(2-\g)\a \ge 2\b+ 2\k\,.
\end{equation}

It remains to bound $\D Z$ in $L^\infty$ and $L^2$. Here there are 2 types of
rounds: ones in which we choose a
triangle that involves $u$ or $v$ and a vertex in $ Y_{u,v}$ (there are $O(n^2 p^4)$ such triangles) and ones where we choose any other triangle, in which case $Y_{u,v}$ is unchanged.  The former event has probability $O( p/n)$ and leads to an $O(1)$ change in $Z$ while the
latter gives a variation in $Z$ of order $O(p/n)$ due to the $-6/n^{2}$ change in $p$. Therefore, $|\D Z| =O(1)$ and $\E[(\D Z)^2\mid\cF_i] = O(p/n)$.

Let $p_0 = p(i_0) = 1 - 6i_0/n^2$. By the definition of $i_0$ and the fact that $|\D Z|=O(1)$,
\[ Z(i_0) \leq (\ha-\a) n^{1/2} (p_0)^{\hg} \Phi + O(1) \leq -\tfrac12 \a n^{1/2} (p_0)^{\hg} \log n \]
(the last inequality holds for large $n$ as the final expression clearly tends to $\infty$ with $n$), and therefore the supermartingale $ S_j = Z\big((i_0 + j) \wedge \tau_P\big) $ satisfies
\begin{align*}
&S_0 \leq -\tfrac12\a n^{1/2} (p_0)^{\hg} \log n \,,\qquad
&\max_j |S_{j+1}-S_j| = O(1)\,, \qquad
&\sum_{j}\E\big[(S_{j+1}-S_j)^2\mid \cF'_j\big] = O \left( n (p_0)^{2} \right)\,,
\end{align*}
where $\cF'_j = \cF_{i_0+j}$. Since $n^{1/2} (p_0)^{2-\hg} \geq c \log^5 n$ for $p_0\geq p^\star$ due to Eq.~\eqref{eq-p-star} we have $|S_0|\max_j |S_{j+1}-S_j| = O\big(n (p_0)^2\big)$ and so Theorem~\ref{thm:concentrate} yields that for some fixed $c>0$
\[ \Pr\left(\cup_j \{S_j \geq 0\}\right) \leq \exp
\left( - c \frac{  n (p_0)^{2\hg} \log^2 n}{ n (p_0)^{2}} \right)
= \exp\left( - c (p_0)^{2\hg-2} \log^2 n \right) \leq e^{-c \log^2 n}\,, \]
which is sufficiently small to show that w.h.p.\ $X(i) < \a n^{1/2} p^\hg \Phi $ for all $u,v$ and $i$.
The same argument handles the analogous symmetric case of the critical interval
$ \big( -\a n^{1/2} p^{\hg}\Phi\,,\, -\ha n^{1/2} p^{\hg}\Phi \big)$ and shows that w.h.p.\ $X(i) > -\a n^{1/2} p^\hg \Phi $ for all $u,v$ and $i$. This concludes the proof of Eq.~\eqref{Yuvbound}.


\subsection{Edges within a neighborhood ($ T_{u}$) }

The number of edges in the subgraph induced by the neighborhood of $u$ can
change in two ways: Either a vertex is removed from $N_u$ (due to selecting a triangle
of the form $uxy$ with $x\in N_u$) thereby decrementing $T_u$ by all edges incident to it in this induced subgraph, or $N_u$ remains unchanged (upon selecting a triangle that does not include $u$)
and yet some of its inner edges are removed.
The former case will be handled by directly summing over $x\in N_u$, noting there
are $Y_{ux}$ triangles of the form $uxy$ whereas the vertex $x$ is incident to $Y_{ux}$
edges counted in $T_u$. The latter case requires a more delicate treatment, similar to the
one used to study $R_{u,v}$ in Section~\ref{sec:R}. Indeed, the naive approach would be
to sum over edges counted by $T_u$, i.e.\ $xy\in E$ with $x,y\in N_u$, as each of these would
decrease $T_u$ by $1$ upon selecting one of the $Y_{x,y}$ triangles incident to it.
However, the cumulative error in this approach (summing the co-degree errors for each edge in $T_u$)
would be quite substantial as it completely ignores the effect of averaging the co-degrees over $T_u$. To take advantage of this point we will use our estimates for the random variables $R_{x,u}$,
which incorporate this averaging effect. Namely,
\begin{equation}\label{eq-Tu-sum}
\E\left[ \D T_u \mid \cF_i\right] = -\frac{1}{Q}\sum_{x\in N_u}Y_{u,x}^2 + \frac{ T_u}{Q}
-\frac{1}{2Q}\sum_{x\in N_u}R_{x,u}\,.
\end{equation}
Here the first two terms accounted for triangles lost due to edges of the form $ux$ (each is chosen with probability $Y_{u,x}/Q$ and eliminates $Y_{u,x}$ triangles from $T_u$, hence the first term, yet in this way each $xy\in T_u$ is double counted, hence the second correcting term). The last term
counted triangles of the form $xyz$ where $x,y,z \in N_u$ as well as ones of the form $xyz$ where $x,y \in N_u$ and $z \not\in N_u \cup \{v\}$. For a given $x$ this corresponds to $R_{x,u}$ (which we recall counts ordered such pairs $(y,z)$, as needed since having $z\in N_{x,u}$ would impact two edges in $T_u$), and the final factor of $\frac12$ makes up for the double count over all $x\in N_u$.

We evaluate the last term in~\eqref{eq-Tu-sum} using the bounds~\eqref{Ruvbound} to get
\begin{align}
-\frac{1}{2Q}\sum_{x\in N_u} R_{x,u} &=
-\frac{1}{2Q} \sum_{x\in N_u}  \left( p Y_{x} Y_{x,u} \pm  \b n^{3/2}p^{2+\hg}\Phi \right)=-\frac{p}{2Q} \sum_{x\in N_u}  \left(Y_x Y_{x,u}\right) \pm  (1+O(\L))\frac{3\b p^{\hg}\Phi}{n^{1/2}}  \,.\label{eq-D-Tu-3-4}
\end{align}

The first sum in~\eqref{eq-Tu-sum} can be estimated by Lemma~\ref{lem:jpatrick} (noting that $\sum_{x\in N_u}Y_{ux}=2T_u$), and so
\begin{align}
\E\left[ \D T_u \mid \cF_i\right] = &-\frac{1}{Q}\bigg(\frac{4T_u^2}{Y_u} \pm (4+O(\L))\a^2n^2p^{1+2\hg}\Phi^2 \bigg)
-\frac{p}{2Q} \sum_{x\in N_u}(Y_x Y_{x,u})
 \pm  (1+O(\L))\frac{3\b p^{\hg}\Phi}{n^{1/2}} \,.\label{eq-D-X-Tu}
\end{align}
As usual
set
\[ X(i)=T_u - pY_u^2/2\]
and consider $\D\left(p Y_{u}^2/2\right)$. Observe that $Y_u$ changes if and only if the triangle selected is of the form $uxy$ with $x,y\in N_{u}$, in which case it decreases by 2. Hence,
$\E[\Delta Y_u \mid \cF_i] =-2T_u/Q$ and $\D (Y_{u}^2) = (2Y_{u}-2)\D Y_u$, and putting these together we get
\begin{align*}
\E\left[\D\left(p Y_{u}^2/2\right)\mid \cF_i\right] &= - \frac6{n^2} \frac{Y_{u}^2}2 - \frac12\bigg(p-\frac6{n^{2}}\bigg)\left(2 Y_{u} -2\right) \frac{2T_u}Q
=- \frac{3Y_{u}^2}{n^2} - \frac{2pT_u Y_u}Q + O\left(\frac{p}{n}\right)\,.
 \end{align*}
Combining this estimate with~\eqref{eq-D-X-Tu} and the bound~\eqref{Yubound} for $Y_u$ gives
\begin{align*}
\E[\D X\mid \cF_i] = &-\frac{4T_u}{QY_u}\left(T_u-p\frac{Y_u^2}{2}\right)-\frac{p}{2Q}\bigg(\sum_{x\in N_u}Y_{x,u}\left(np\pm \k n^{1/2}p^{\hg-1}\Phi \right)\bigg)+
\frac{3 Y_u^2}{n^2}\\
& \pm  \left( 1 + O( \L ) \right) \left[ \frac{3\b p^\hg \Phi }{ n^{1/2}} + \frac{ 24 \a^2 p^{2\hg -2} \Phi^2}{n}
\right] 
\end{align*}
By Eq.~\eqref{eq-p-gamma-comparison} we have $n^{-1} p^{2\hg-2}\Phi^2 = o(\L n^{-1/2} p^{\hg}\Phi)$ for all $p \geq p^\star$, thus the above expression involving $\a^2$ can be absorbed into the $O(\L)$ error-term of the expression involving $\b$. Furthermore, since $4T_u / (Q Y_u) = (12+O(\L))/(n^2p)$ and
\[ \frac{p}{2Q}\bigg(\sum_{x\in N_u}Y_{x,u}\left(np\pm \k n^{1/2}p^{\hg-1}\Phi \right)\bigg)-
\frac{3 Y_u^2}{n^2} =
\frac{6}{n^2 p} \left(T_u -p Y_u^2/2\right) \pm (3+O(\L))\frac{\k p^{\hg}\Phi}{n^{1/2}}  + O\left( \frac{p^{2\hg-2}\Phi^2}n\right)
\]
we can conclude that
\[
\E[\D X\mid \cF_i] = -\frac{18+O( \L )}{n^2p}X
\pm  \left( 1 + O( \L ) \right) \frac{3(\b+\k) p^\hg \Phi }{ n^{1/2}} \,.
\]
Now we consider the upper critical interval for $ T_u$ given by
\[
I_T = \left( \hm n^{3/2} p^{1+\hg}\Phi~,~ \m n^{3/2} p^{1+\hg}\Phi \right) \qquad \mbox{ where } \qquad \hm = \left( 1- \log^{-1} n \right) \m\,,
\]
and as before let $i_0$ be the first round in which $X(i) \geq \hm n^{3/2}p^{1+\hg}\Phi$, define the stopping time $\tau_T = \min\{ i \geq i_0 : X(i) < \hm n^{3/2}p^{1+\hg}\Phi\}$ and consider
\[ Z(i) = X(i) - \m n^{3/2} p^{1+\hg} \Phi \,.\]
Exactly the same argument of~\eqref{eq-D-2nd-Z-term} gives
\[ \D\left[\m n^{3/2} p^{1+\hg} \Phi\right] = \left(6+O(n^{-2})\right)\left(p^{1+\hg} - (1+\hg) p^{\hg}\right)\m n^{-1/2} \Phi \]
and together with the above upper bound on $\E\left[ \D X\mid \cF_i\right]$ we get
\begin{align*}
\E\left[ \D Z\mid \cF_i\,,\,\tau_T>i\right]
& \leq \big[-6\hm + \b + \k + (2+2\hg)\m +O(\L)\big]\frac{ 3p^{\hg}\Phi }{ n^{1/2}}
+ \big[- \m + O(\L)\big] \frac{ 6p^{1+\hg} \Phi }{ n^{1/2}} \\
&\leq \big[-(4-2\g)\m + \b + \k \big]\frac{ 3p^{\hg}\Phi }{ n^{1/2}}
+ \big( - \m + O(\L)\big) \frac{ 6p^{1+\hg} \Phi }{ n^{1/2}} \,,
\end{align*}
where we used the fact that $\g - \hg = 6/\log n$ to absorb both $4(\m-\hm)=4\mu/\log n$ and the $O(\L)$-term for large $n$. In particular, $S_j = Z((i_0+j) \wedge \tau_Y)$ is indeed a supermartingale provided that
$\Lambda=o(\mu)$ and
\begin{equation}
\label{eq:keyT}
(4-2\g) \m \ge \b+\k\,.
\end{equation}

Having established an $L^1$ bound on $\D S$ it remains to consider the corresponding $L^2,L^\infty$ bounds.
If we choose a triangle of the form $uxy$, an event that has probability $O(1/n)$,
then $T_u$ decreases by $Y_{ux}+Y_{uy}-1$ while $Y_u$ decreases by $2$, hence the change in $Z$ in this case is
at most the $ O( n^{1/2} p^\hg \log^2n) $ due to the error-terms in our approximation for the degrees and co-degrees.
The probability that we choose a triangle that does not contain $u$ yet includes an edge in $T_u$ is
$O( p^2)$ and selecting such a triangle changes $Z$ by $O(1)$.  The choice of
any other triangle changes $Z$ by $O(p^2)$ due to the change in $p$.
Altogether, $|\D Z| = O\big(n^{1/2}p^\hg \log^2 n\big)$ and $\E[(\D Z)^2\mid \cF_i] = O(p^{2\hg} \log^4 n)$.
Let $p_0 = p(i_0) = 1 - 6i_0/n^2$ and recall that the definition of $i_0$ and our bound on $|\D Z|$ ensure that
\[ Z(i_0) \leq (\hm-\m) n^{3/2} (p_0)^{1 + \hg} \Phi  + O\big( n^{1/2} (p_0)^{\hg} \log^2n \big) \leq -\tfrac12 \m n^{3/2} (p_0)^{1 + \hg} \log n \,,\]
where the factor of $\frac12$ absorbs the $O\big( n^{1/2} (p_0)^{\hg} \log^2 n\big)$-term since $n p_0 \log^{-1} n \geq n^{1-o(1)}p^\star > \sqrt{n}$ for large enough $n$. It then follows the supermartingale $S_j=Z((i_0+j)\wedge\tau_T)$ satisfies
\begin{align*}
&S_0 \leq -\tfrac12 \m n^{3/2} (p_0)^{1+\hg} \log n\,, \qquad\qquad
\max_j |S_{j+1}-S_j| = O\big(n^{1/2}p^\hg \log^2 n\big)\,,\\
&\sum_{j}\E\big[(S_{j+1}-S_j)^2\mid \cF'_j\big] = O \big( n^{2} (p_0)^{1+2\hg}\log^4 n \big)\,,
\end{align*}
where $\cF'_j = \cF_{i_0+j}$. Here $|S_0|\max_j |S_{j+1}-S_j| = o\big(n^{2} (p_0)^{1+2\hg}\log^4 n\big)$ due to one extra log factor between these expressions and therefore Theorem~\ref{thm:concentrate} establishes that for some fixed $c>0$
\[ \Pr\left(\cup_j \{S_j \geq 0\}\right) \leq \exp
\left( - c \frac{  n^3 (p_0)^{2+2\hg} \log^2 n}{ n^2 (p_0)^{1+2\hg}\log^4 n} \right)
= \exp\left( - c n p_0 \log^{-2} n \right) \leq e^{- \sqrt{n}}\,. \]
We conclude that w.h.p.\ $X(i) < \m n^{3/2} p^{1+\hg} \Phi $ for all $u$ and $i$, and the same argument shows that w.h.p.\
$X(i) > -\m n^{3/2} p^\hg \Phi $ for all $u$ and $i$. This concludes the proof of Eq.~\eqref{Tubound}.


\subsection{Vertex degrees ($ Y_{u}$) }
The analysis of the degrees will be straightforward using our estimate~\eqref{Tubound} for $T_u$, the number of inner edges in the neighborhood of a vertex $u$, since $Y_u$ changes iff the triangle selected is of the form $uxy$ (in which case it decreases by 2). Indeed, setting
\[ X(i)=Y_u-np\,,\]
our bounds on $T_u$ and $Q$ imply that
\begin{align*}
\E[\Delta X\mid \cF_i] & =-\frac{2T_u}{Q}+\frac{6}{n} = -\frac{p Y_u^2 \pm 2\m n^{3/2}p^{1+\hg}\Phi}{\frac16n^3p^3\pm\a^2n^2p^{2\hg-1}\Phi^2}+\frac{6}{n}\\
& = -\frac{ p\left( X^2+(np)^2 + 2npX \right)}{ \frac16 n^3p^3}+\frac{6}{n} \pm (12+O(\L))\frac{ \m p^{\hg-2}\Phi}{ n^{3/2}}
+O \left( \frac{ p^{2\hg -4}\Phi^2}{ n^2} \right) \\
&= - \frac{ 12 }{ n^2 p}X \pm ( 12 + O(\L))\frac{ \mu p^{\hg - 2}\Phi}{ n^{3/2}}\,,
\end{align*}
where the term $pX^2 / (n^3p^3)$ was absorbed into the $O(p^{2\hg-4}\Phi^2 / n^2)$-term since $|X|\leq \k \sqrt{n}p^{\hg-1}\Phi$ by~\eqref{Yubound}, and this latter error-term was thereafter omitted as it is $o(\L p^{\hg-2}\Phi/n^{3/2})$ by~\eqref{eq-p-gamma-comparison}.

Now consider the upper critical interval for $Y_u$ given by
\[
I_Y' = \left( \hk n^{1/2} p^{\hg-1}\Phi~,~ \k n^{1/2} p^{\hg-1}\Phi \right) \qquad \mbox{ where } \qquad \hk = \left( 1- \log^{-1} n \right) \k\,,
\]
and as before assume that $i_0$ is the first time at which $X(i_0) \geq \hk  n^{1/2} p^{\hg-1} \Phi$ and define the stopping time $\tau'_Y = \min\{i>i_0 : X(i) < \hk  n^{1/2} p^{\hg-1} \Phi\}$. With these definitions, the variation in the variable
\[ Z(i)=X(i)-\k n^{1/2}p^{\hg-1} \Phi\]
consists of $\D X$ as well as $ \D\left[-\k n^{1/2} p^{\hg-1} \Phi\right] \leq \left(6+O(n^{-2})\right)(\hg-1) \k p^{\hg-2} n^{-3/2} \Phi $, hence
\[\E\left[ \D Z\mid \cF_i\,,\,\tau_Y'>i\right] \leq
\big[-2\hk+2\m +(\hg-1)\k + O(\L)\big]\frac{6p^{\hg-2}\Phi}{n^{3/2}}  \leq \big[-(3-\g)\k+2\m \big]\frac{6p^{\hg-2}\Phi}{n^{3/2}}\,,
\]
where we used the fact that $\g-\hg = 6/\log n$ eliminates the term $2(\k-\hk)=2\k/\log n$ and $O(\L)$-term for large $n$. Hence, $S_j = Z((i_0+j) \wedge \tau_Y')$ is a supermartingale as long as
\begin{equation}
\label{eq:keyY}
(3-\g) \k \ge 2\m\,.
\end{equation}
Furthermore, in each round we either select a triangle incident to $u$, an event which has probability $O(1/n)$ and changes $Z$ by $O(1)$,
or we do not affect $Y_u$ and thus change $Z$ by $O(1/n)$ due to the change in $p$. Thus, $|\D Z| = O(1)$ while $\E[(\D Z)^2\mid\cF_i] = O(1/n)$, and we conclude that for large enough $n$ the supermartingale $S_j$ has the following attributes:
\begin{align*}
&S_0 \leq -\tfrac12 \k n^{1/2} (p_0)^{\hg-1} \log n ~,~\quad \max_j |S_{j+1}-S_j| = O(1)~,~\quad \sum_{j}\E\big[(S_{j+1}-S_j)^2\mid \cF'_j\big] = O ( n p_0 )\,,
\end{align*}
where $p_0 = p(i_0) = 1 - 6i_0/n^2$ and $\cF'_j = \cF_{i_0+j}$. Since $\sqrt{n}(p_0)^{2-\hg} \geq \sqrt{n}(p^\star)^{2-\hg} \geq c \log^5 n$ we deduce that $|S_0|\max_j |S_{j+1}-S_j| = o(n p_0)$ and thus Theorem~\ref{thm:concentrate} establishes that for some fixed $c>0$
\[ \Pr\left(\cup_j \{S_j \geq 0\}\right) \leq \exp
\left( - c \frac{  n (p_0)^{2\hg-2} \log^2 n}{ n p_0} \right)
= \exp\left( - c p_0^{2\hg-3} \log^{2} n \right) \leq e^{- c\log^2 n}\,. \]
Altogether, w.h.p.\ $X(i) < \k n^{1/2} p^{\hg-1} \Phi $ for all $u$ and $i$, and similarly
$X(i) > -\k n^{1/2} p^{\hg-1} \Phi $ w.h.p.\ for all $u$ and $i$. This completes the proof of Eq.~\eqref{Yubound}.


\subsection{Co-degree of triples ($ Y_{u,v,w}$) }
\label{sec:Yuvw}
We will prove the following result from which~\eqref{Yuvwbound} will readily follow:
\begin{equation}
  \label{eq:Yuvwbound-str}
  Y_{u,v,w} = n p^3 \pm \sqrt{n p^3 \log^5 n}\quad \mbox{for all $u,v,w$ and $p \geq p_1 := n^{-1/3} \log^{5/3} n$}\,.
\end{equation}
%
Define
\[X(i)=Y_{u,v,w}(i)-np^3.\]
We have
\begin{align*}
 \E[\D X\mid \cF_i]
& = - \sum_{ x \in N_{u,v,w}} \frac{ Y_{ x,u} + Y_{x,v} + Y_{x,w} -\one_{\{uv\in E\}}-\one_{\{uw\in E\}}-\one_{\{vw\in E\}}}{Q}
+ \frac{18 p^2}{n}\\
& \le - \frac{ (X + n p^3)  (3np^2 - \alpha n^{1/2} p^{\hg} \Phi)}{ n^3p^3/6 }  + O \left( \frac{\a^2 p^{2\hg -2} \Phi^2}{ n^2} \right)  + \frac{ 18 p^2}{n}  \\
& = - \frac{ 18}{ n^2 p}X + \frac{ 6\a p^{\hg}\Phi}{ n^{3/2}} + O \left( \frac{\a p^{\hg -3/2} \Phi \log^2 n}{ n^2} \right)
\end{align*}
Now suppose that $i_0$ is a first round at which $ X(i_0) >  \tfrac23 \sqrt{ n p^3 \log^5 n} $, that is $X$ enters the
critical interval \[I_Y'' = \Big(\; \tfrac23 \sqrt{ n p^3 \log^5 n}~,~  \sqrt{ n p^3 \log^5 n}\;\Big)\,.\]
Setting $\tau_Y'' = \min\{i>i_0: X(i) < \tfrac23\sqrt{np^3\log^5n}\}$ and
\[ Z(i) = X(i) - \sqrt{ n p^3 \log^5 n}\]
we get
\begin{align*}
\E\left[ \D Z \mid \cF_i\,,\,\tau''_Y>i\right]  &\le  -\frac{12 p^{1/2} \log^{5/2} n}{n^{3/2}} +\frac{ 6\a p^{\hg}\Phi}{ n^{3/2}}
+ O \bigg( \frac{\a p^{\hg -3/2} \Phi \log^2 n}{ n^2} \bigg) + \frac{ 9 p^{1/2} \log^{5/2} n}{ n^{3/2}} \\
&= \frac{6\a p^\hg \Phi- 3p^{1/2}\log^{5/2} n  }{n^{3/2}} + O \bigg( \frac{\a p^{\hg -3/2} \Phi \log^2 n}{ n^2} \bigg) \,.
\end{align*}
Since $\hg = \frac12 - O(1/\log n)$ it follows that
 $p^\hg \Phi = O(p^{1/2} \log^2 n) = o( p^{1/2} \log^{5/2} n)$ and hence the first term
in the above r.h.s.\ is equal to $-(3+o(1)) \big(p \frac{\log^5 n}{n^3}\big)^{1/2}$. As for the second term there, recall that $p \geq p^\star$ and so by~\eqref{eq-p-star} we have $p^{4-2\hg} \geq 6\alpha^2 \Phi^2 n^{-1} \log^{6}n$. In particular,
\[ \frac{p^{\hg-3/2}\Phi \log^2 n}{n^2} = \bigg(\frac{p\log^5 n}{n^3}\bigg)^{1/2} \frac{\Phi}{n^{1/2} p^{2-\hg}\log^{1/2} n} \leq \bigg(\frac{p\log^5 n}{n^3}\bigg)^{1/2} \frac{1}{\sqrt{6}\alpha \log^{7/2} n}\,,
 \]
and altogether we conclude that
\[ \E\left[ \D Z \mid \cF_i\,,\,\tau''_Y>i\right]  \leq -(3+o(1)) \frac{p^{1/2}\log^{5/2} n}{n^{3/2}} < 0\,,\]
where the last inequality holds any for sufficiently large $n$ and confirms that $S_j=Z((i_0+j)\wedge\tau_Y'')$ is a supermartingale.
Moreover, the equation that specified $\E[\D X\mid \cF_i]$ shows that if the triangle selected goes through $u,v$ or $w$ and a vertex in $ N_{u,v,w}$,
which happens with probability $O(p^2/n)$, then
the change in $Z$ is $O(1)$, and otherwise the change in $Z$ is $O(p^2/n)$.
Hence, $|\D Z| = O(1)$ while $\E[(\D Z)^2\mid\cF_i] = O(p^2/n)$, thus the supermartingale $S_j$ has the following attributes:
\begin{align*}
&S_0 \leq -\tfrac14 \sqrt{ n (p_0)^3 \log^5 n}~,~\quad \max_j |S_{j+1}-S_j| = O(1)~,~\quad \sum_{j}\E\big[(S_{j+1}-S_j)^2\mid \cF'_j\big] = O \big( n (p_0)^3 \big)\,,
\end{align*}
where the factor $\frac14$ in the first expression (as opposed to $\frac13$) treated the potential $O(1)$ deviation of $Z(i_0)$ from the lower endpoint of the critical interval.
Noting that $p_0 \geq p_1 = \big(\frac{\log^{5}n}{n}\big)^{1/3}$ and hence $|S_0|\max_j |S_{j+1}-S_j| = O\big(n(p_0)^3\big)$, an application of Theorem~\ref{thm:concentrate} yields that for some fixed $c>0$,
\[ \Pr\left(\cup_j \{S_j \geq 0\}\right) \leq
\exp \left( - c \frac{ n (p_0)^{3} \log^5 n}{ n (p_0)^3 } \right) = \exp\left(-c \log^5 n \right)\,. \]
By the usual union bound over vertices and rounds we now conclude that w.h.p.\ $X(i) < \sqrt{np^3\log^5 n}$ for all $u,v,w$ and $i$, and similarly
$X(i) > -\sqrt{np^3\log^5 n} $ w.h.p.\ for all $u,v,w$ and $i$, thus completing the proof of Eq.~\eqref{Yuvwbound}.


\subsection{Number of triangles ($Q$)}
\label{sec:Q}

In~\cite{BFL} it was shown (see Theorem~2 there) that
$ Q(i)\leq \frac{n^3p^3}{6}+\frac{1}{3}n^2p $
throughout the process w.h.p., hence it only remains to prove the lower bound in~\eqref{Qbound}.
Let
\[ X(i) = Q - \frac{n^3p^3}{6}\]
and recall that due to~\eqref{eq-Q-one-step} we have
\[ \E[\D Q \mid \cF_i] = - \sum_{ xyz \in Q} \frac{ Y_{x,y} + Y_{x,z} + Y_{y,z} - 2}{Q}
\geq - \frac{1}{Q} \sum_{ xy \in E} Y_{x,y}^2\,.\]
To bound $\E[\D X \mid \cF_i]$ from below we will thus need an upper bound on $\sum_{ xy \in E} Y_{x,y}^2$.
Recall that $\sum_{x y\in E} Y_{x,y} = 3Q$ and that $Y_{x,y} = np^2 \pm \a n^{1/2}p^{\hg} \Phi$ by Eq.~\eqref{Yuvbound}, hence we can apply Lemma~\ref{lem:jpatrick} together with the fact that $|E(i)| = n^2 p/2 - n/2$ to obtain that
\begin{align*}
\frac1Q \sum_{xy\in E} Y_{x,y}^2 &\leq \frac1{Q}\bigg(\frac{9Q^2}{|E|} + 4|E| \a^2 n p^{2\hg}\Phi^2\bigg) \leq \frac{18 Q}{n^2p\big(1-\frac1{np}\big)}+
\frac{2n^2 p}{\big(\frac16+O(\L)\big)n^3p^3} \a^2 n p^{2\hg}\Phi^2 \\
&\leq  \frac{18 Q}{n^2p}+ O(p)+ (12+O(\L))\a^2 p^{2\hg-2}\Phi^2
= \frac{18 Q}{n^2p}+ (12+O(\L))\a^2 p^{2\hg-2}\Phi^2
\,,
 \end{align*}
 where in the last equality we absorbed the $O(p)$-term into the $O(\L)$ error-term factor of the last expression since the facts $\hg\leq 1$ and $\L \Phi^2\geq c\log^3n$ imply that $p^{2-2\hg} = o(\L \Phi^2)$, i.e.\ $\L p^{2\hg-2}\Phi^2\to\infty$. Adding this to our estimate for $\E[\D Q \mid \cF_i]$ while observing that $\D (-\frac16 n^3p^3) = 3np^2 + O(p/n)$ yields
 \[ \E[\D X \mid \cF_i] \geq - \frac{18 Q}{n^2p} - (12+O(\L))\a^2 p^{2\hg-2}\Phi^2 + 3np^2 + O(p/n) \,.\]
As before we incorporate the $O(p/n)$ term into the $O(\L)$ error and using the definition of $X$ we can then rewrite the above as an upper bound on $\D(-X)$, as follows:
 \[ \E[\D (-X) \mid \cF_i] \leq \frac{18 }{n^2p}X + (12+O(\L))\a^2 p^{2\hg-2}\Phi^2 \,.\]
Now assume that $i_0$ is the first round where $X$ drops below $ -\ha^2 n^{2} p^{2\hg-1}\Phi^2$, i.e.\ enters the interval
\[
I_Q = \left( -\a^2 n^{2} p^{2\hg-1}\Phi^2~,~ -\ha^2 n^{2} p^{2\hg-1}\Phi^2 \right) \qquad \mbox{ where } \qquad \ha = \left( 1- \log^{-1} n \right)^{1/2} \a\,.
\]
Further let $\tau_Q = \min\{i>i_0: X(i) > -\ha^2 n^{2} p^{2\hg-1}\Phi^2 \}$ and
\[ Z(i) = -X(i) - \a^2 n^2 p^{2\hg-1}\Phi^2 \,.\]
Since $ \D\left[-\a^2 n^{2} p^{2\hg-1} \Phi^2\right] \leq \left(6+O(n^{-2})\right)(2\hg-1) \a^2 p^{2\hg-2} \Phi^2$, the upper bound on $\D(-X)$ gives
\[\E\left[ \D Z\mid \cF_i\,,\,\tau_Q>i\right] \leq
\big[-3\ha^2 +2\a^2 + (2\hg-1)\a^2 + O(\L)\big]6p^{2\hg-2}\Phi^2  \leq -12(1-\g)\a^2 p^{2\hg-2}\Phi^2\,,
\]
where the last inequality used the term $\g-\hg = 6/\log n$ to both cancel the $O(\L)$-term and replace $3\ha^2$ by $3\a^2$. As $\g<1$ we deduce that $S_j = Z((i_0+j) \wedge \tau_Q)$ is indeed a supermartingale.

Next consider the one-step variation of $Z$. Denoting the selected triangle in a given round by $xyz$, the change in $Q$ following this round is at most $Y_{x,y}+Y_{x,z}+Y_{y,z}$ and in light of our co-degree estimate~\eqref{Yuvbound} this expression deviates from its expected value of $3np^2$ by at most $3\a n^{1/2}p^{\hg} \log^2 n$. In particular, $|\D Z| = O\big( \sqrt{n}(p_0)^{\hg} \log^2 n \big)$ and letting $p_0 = p(i_0) = 1 - 6i_0/n^2$ this ensures that
\[ Z(i_0) \leq (\ha^2-\a^2) n^{2} (p_0)^{2\hg-1} \Phi^2  + O\big( \sqrt{n}(p_0)^{\hg} \log^2 n \big) \leq -\tfrac12 \a^2 n^{2} (p_0)^{2\hg-1} \log^3 n \,,\]
where the last inequality holds for large $n$ since $n^{3/2} (p_0)^{\hg-1}\log n $ tends to $\infty$ with $n$. With at most $n^2 p_0$ steps remaining until the process terminates, Hoeffding's inequality establishes that for some fixed $c>0$,
\[ \Pr\left(\cup_j \{S_j \geq 0\}\right) \leq \exp
\left( - c \frac{  (n^2 (p_0)^{2\hg-1} \log^3 n)^2}{ n^2 p_0 (\sqrt{n} (p_0)^{\hg}\log^2 n)^2} \right)
= \exp\left( - c n p_0^{2\hg-3} \log^{2} n \right) \leq e^{- n}\,. \]
We conclude that w.h.p.\ $X(i) > -\a^2 n^2 p^{2\hg-1} \Phi^2$  for all $i$, completing the proof of Eq.~\eqref{Qbound}.

This completes the proof of Theorem~\ref{thm:enchilada}.  \qed

\end{document}